\documentclass[11pt]{article}
\usepackage{pstricks,pst-node}
\usepackage{amsmath}
\usepackage{amssymb}
\usepackage{amsthm}
\usepackage{amsfonts}
\usepackage{graphicx,fancyheadings}
\newcommand{\al}{\alpha}
\newcommand{\be}{\beta}

\newcommand{\ep}{\epsilon}

\newcommand{\la}{\lambda}

\newcommand{\qf}{Q_F}
\newcommand{\lt}{\rightsquigarrow}

\newcommand{\lton}{\longrightarrow^{\hspace{-4mm}n}\hspace{2mm}}

\newcommand{\yn}{Y_n}
\newcommand{\xin}{X_{i,n}}
\newcommand{\xn}{X_n}
\newcommand{\tmn}{T_{m,n}}

\newcommand{\lanm}{\lambda_{n,m}}
\newcommand{\Gin}{\Gamma_{i,n}}
\newcommand{\Gn}{\Gamma_n}
\newcommand{\muin}{\mu_{i,n}}

\newgray{gray9}{.9}
\newgray{gray8}{.8}
\newgray{gray85}{.85}

\newtheorem{thm}{Theorem}
\newtheorem{lem}[thm]{Lemma}
\newtheorem{cor}[thm]{Corollary}

\title{Splitting pairs and the number of clusters generated by random pair incompatibilities}

\author{Damien Pitman 
\smallskip
\thanks{The author would like to thank Janko Gravner for suggesting this model and for his help in revising the paper. Also, this research was funded in part by the National Science Foundation grant DMS-0204376.}\\
Mathematics Department, University of California, Davis\\
pitman@math.ucdavis.edu}

\begin{document}

\maketitle

\thispagestyle{empty}

\begin{abstract}

We consider a random fitness landscape on the space of haploid diallelic genotypes with $n$ genetic loci, where each genotype is considered either inviable or viable depending on whether or not there are any incompatibilities among its allele pairs. We suppose that each allele pair in the set of all possible allele pairs on the $n$ loci is independently incompatible with probability $p=c/(2n)$. We examine the connectivity of the viable genotypes under single locus mutations and show that, for $0<c<1$, the number of clusters of viable genotypes in this landscape converges weakly (in $n$) to $N=2^\Psi$ where $\Psi$ is Poisson distributed; while for $c>1$, there are no viable genotypes with probability converging to one. The genotype space is equivalent to the hypercube $Q^n$ and the viable genotypes are solutions to a random 2-SAT problem, so the same result holds for the connectivity of solutions in $Q^n$ to a random 2-SAT problem.

\end{abstract}
\begin{section}{Introduction}

The space of diallelic haploid genotypes on $n$ genetic loci is in bijection with the set of binary strings of length $n$ or the vertices of the $n$-cube $Q^n$. Let $[n]=\{1,\ldots,n\}$ be the set of genetic loci. At each locus $i$, let the two alleles be denoted $0_i$ and $1_i$. We define a pair of alleles $(x,y)$ on distinct loci to be {\em incompatible} if the existence of these alleles in a genotype is lethal. Suppose we are given a list of incompatibilities $L$. We say that a genotype is {\em viable} if none of its allele pairs is on the list $L$. 

This model is an example of a {\em fitness landscape}. The notion of fitness landscapes was introduced by a theoretical evolutionary biologist, Sewall
Wright in \cite{Wri32} (see also \cite{Kau93,Gav}). The study of fitness landscapes has proved extremely useful both in
biology and well outside of it. In the standard interpretation, a fitness landscape is a relationship
between a set of genes (or a set of quantitative characters) and a measure of fitness
(e.g. viability, fertility, or mating success). In Wright's original formulation the set of 
genes (or quantitative characters) is the property of an individual. However, the notion of 
fitness landscapes can be generalized to the level of a mating pair, or even a population of
individuals. For a comprehensive introduction to fitness landscapes see~\cite{Gav}; or for descriptions more closely aligned with the ideas here, see~\cite{GPG} or~\cite{GG}.

This model is easily translated into propositional logic. To each locus $1 \leq i \leq n$ we associate a boolean variable $\xi_i$. The two truth values of $\xi_i$ are referred to as the {\em positive literal} and the {\em negative literal} of $\xi_i$. We make the convention that the allele $0_i$ is the negative literal and $1_i$ is the positive literal. To identify a literal in this model we need both locus and value, but we put the locus in the subscript. 
The negation of any literal $x$ is denoted $\overline{x}$. Thus, $\overline{0_i}=1_i$ and $\overline{1_i}=0_i$. Alleles at distinct loci or literals of distinct variables are said to be {\em strictly distinct}. An assignment of $n$ strictly distinct literals to $n$ boolean variables is referred to as a truth assignment to the $n$ variables. Thus, each genotype is equivalent to a truth assignment.

If $(x,y)$ is an incompatibility, then a genotype that has allele $x$ must also have allele $\bar{y}$ to be viable, and a genotype with allele $y$ must have $\overline{x}$. Thus, there are two implications naturally associated to the incompatibility $(x,y)$; $x \Rightarrow \bar{y}$ and $y \Rightarrow \overline{x}$. It follows that the incompatibility $(x,y)$ is equivalent to the 2-clause or disjunction $\overline{x} \vee \bar{y}$. Let $c_1,\ldots,c_m$ be 2-clauses with literals chosen from the $n$ variables. Then 
$$F=c_1 \wedge \ldots \wedge c_m $$
is a {\em 2-formula}. If there is some literal assignment to the $n$ variables for which the formula $F$ is TRUE, then it is said to be {\em satisfied} or SAT. Determining if a formula is satisfiable is known as the 2-SAT problem.

Given the equivalence of 2-clauses and incompatibilities, the set of lists of $m$ incompatibilities and the set of 2-formulae of length $m$ are in bijection. Moreover, if $L$ is a list of incompatibilities and $F$ the corresponding 2-formula, then a genotype having no allele pairs on $L$ is equivalent to the associated truth assignment satisfying the 2-formula $F$. Thus, every viable genotype determined by $L$ is a solution to the 2-SAT problem corresponding to $F$ and vice versa. For the remainder of this paper, we let $F$ refer to either the formula or the list of incompatibilities without distinction.

We are interested in the connectivity of viable genotypes through single locus mutations. If two viable genotypes are connected by a path of single locus mutations we say that they are in the same cluster. Thus, each cluster of viable genotypes corresponds to a collection of genotypes which can each evolve into any other genotype in the cluster without passing through an inviable genotype and without altering two or more alleles simultaneously. We let the edges of $Q^n$ represent single locus mutations, so the clusters of viable genotypes determined by a formula $F$ correspond to maximal connected components in the subgraph, $Q_F = Q_F^n \subset Q^n$.

A random subgraph can be chosen by supposing that each pair of alleles is incompatible independently with probability $p$. We denote the resulting set of incompatibilities by $F_p=F_p^n$ and let $Q_{F_p}=Q_{F_p}^n$ be the random subgraph that is induced by the viable genotypes. Since it will be clear from context whether or not the set of incompatibilities is randomly determined, we simply denote this set by $F$ (or $F^n$).

The main result of this paper concerns the random variable $N_n$ which counts the number of clusters of viable genotypes in the random subgraph $Q_F$. 
\begin{thm}\label{thm:noocl}
Let $F$, $Q_F$ and $N_n$ be randomly determined as above. If $p=c/(2n)$ and $c < 1$, then $N_n$ converges weakly to $N=2^\Psi$ where $\Psi$ is a Poisson random variable with mean 
$$\la= - \frac{1}{2} \left( \ln{(1-c)}+c \right) .$$
In particular, the probability that there is a unique cluster converges to 
$$e^{-\la}=\sqrt{(1-c)e^c}.$$
\end{thm}

Determining the probability that there are any viable genotypes in this model is the random 2-satisfiability or random 2-SAT problem where each possible 2-clause is included with probability $p$. The results on random 2-SAT found in \cite{CR}~\cite{Go1}~\cite{dlV1}~\cite{BBCKW}~\cite{BD}~\cite{GPG} imply the following statements with probability approaching one as $n$ approaches infinity. 

Let $p=c/(2n)$, where $c>0$ is constant. Then, if $c>1$, $F$ allows no viable genotypes; and if $c<1$, then the number of viable genotypes in $Q_F$ increases exponentially with $n$. 

We should also point out that the theorem here may be related to the clustering results concerning solutions to random 2-SAT found in~\cite{BMW} utilizing the {\em replica method}.

\end{section}
\begin{section}{The digraph}

Every formula $F$ can be associated to a digraph $D_F$ as defined by Aspvall et al \cite{APT}. The vertex set of $D_F$ is the set of alleles 
$$\{ 0_1,\ldots,0_n,1_1,\ldots,1_n \}.$$ 
For each incompatibility $(x,y) \in F$ the directed edges $x \to \bar{y}$ and $y \to \overline{x}$ are included in $D_F$, one edge for each implication. If there is a directed path from $x$ to $\bar{y}$ in $D_F$ it is denoted by $x\lt \bar{y}$ and this path is in $D_F$ if and only if the path $y \lt \overline{x}$ is in $D_F$. 

The set of alleles and edges that can be reached from $x$ is referred to as the out-graph $D^+_F(x)$ and the allele set of $D^+_F$ is denoted
$$L^+(x)=L^+_F(x)=\{y:x\lt y\}.$$
Since edges are equivalent to implications, these are the alleles that must be present in any viable genotype with allele $x$. Of course, a genotype may have all the alleles in $L^+(x)$ and yet still be inviable due to an unrelated incompatibility. 

We are also interested in those alleles that require $x$ for their own existence in a viable genotype; $D^-_F(x)$ denotes the in-graph of $x$. The allele set of the in-graph of $x$ is denoted
$$L^-(x)=L^-_F(x)=\{y:y\lt x\}.$$

Suppose that 
$y \in L^+_F(x) \cap L^-_F(x) $
, so $x \lt y$ and $y \lt x$. In this case $x$ and $y$ are said to be {\em strongly connected}. The set of alleles 
$$C_F(x)=\{y: x \lt y \lt x\}$$
that are strongly connected to $x$ is referred to as the {\em strong component} of $x$. We follow the convention that $x\lt x$ for every vertex $x\in D_F$. Consequently, every $x$ is a member of a strong component and $C_F(x)$ is well defined for every allele. Thus, the set of all $2n$ alleles can be partitioned into strong components of $D_F$ and a partial order can be defined on the set of strong components. Let 
$$C_F(x) \leq C_F(y) \quad \text{if} \quad x \lt y,$$
which also implies that $x' \lt y'$ for any $x' \in L^-_F(x) \supseteq C_F(x)$ and $y' \in L^+_F(y) \supseteq C_F(y)$. 

We think of the alleles of a strong component as depending upon each other for survival; a viable genotype either has all or none of the alleles within a strong component. If a strong component is of size one, i.e., $C_F(x)=\{x\}$, then this is an empty statement and we say that the strong component of $x$ is trivial. Nontrivial strong components represent genetic groupings and the clustering behavior of $Q_F$ is determined by these groupings, so we are quite interested in strong components. However, combinatorial problems relating to strong components can become unwieldy. We will avoid these difficulties by focusing our attention on directed cycles rather than strong components. We will see that in this random setting, strong components tend to be cycles.

Notice that for any strong component $C_F(x)$ there is a unique subgraph of $D_F$ that contains the alleles of $C_F(x)$ along with all the edges between them that are present in $D_F$. Thus, there should be no confusion if we refer to these subgraphs as strong components as well. With this graphical notion of a strong component in mind, we notice that if $x$ and $y$ are distinct alleles in the same strong component of $D_F$, then there is a closed directed walk in $D_F$ containing $x$ and $y$. When a directed walk of length $i$ is allele disjoint, i.e., has $i$ edges and $i$ alleles, then we say that it is a {\em simple} cycle of length $i$, or more succinctly, an $i$-cycle. We refer to a closed walk that contains one or more repeated alleles as a {\em compound cycle}.

Next, for any set of alleles $M$, we define
$$\overline{M}=\{\bar{y}: y \in M \}.$$
When $ M \cap \overline{M} = \emptyset$, we say that $M$ and $\overline{M}$ are {\em complementary}.
Notice that if $y\in C_F(x)$, then 
$\overline{y}\in C_F(\overline{x})$ so $C_F(\overline{x})=\overline{C_F(x)}$ for every allele $x$. Moreover,
$$ C_F(x) \cap \overline{C_F(x)} \neq \emptyset \iff C_F(x)=\overline{C_F(x)} .$$
If there is an allele $x$ such that $C_F(x)=\overline{C_F(x)}$, then there are alleles $y,z \notin \{x,\overline{x}\}$ (not necessarily distinct) such that 
$$x \lt y \lt \overline{x} \lt z \lt x$$
and we say that $x$ (and each $y \in C_F(x)$) is on a {\it contradictory cycle}. This cycle may be compound or simple. Since the allele $x$ depends upon its complement $\overline{x}$ for survival and vice versa, neither allele can be assigned at the corresponding locus and there are no viable genotypes. On the other hand, it can be shown that if there are no contradictory cycles, then there is at least one viable genotype (see \cite{BBCKW} for example). Thus, a formula $F$ allows viable genotypes in $Q_F$ if and only if there is no contradictory cycle in the associated digraph $D_F$. 

\end{section}
\begin{section}{The number of clusters}
\begin{subsection}{A pair of deterministic results}

In this subsection we consider the relationship between a fixed formula $F$ and the structure of $Q_F$. 
Whenever a formula $F$ is satisfiable, there is no allele $x$ such that $C_F(x)=\overline{C_F(x)}$, which implies that for a digraph corresponding to a viable genotype space, strong components come in complementary pairs $(C_F(x),\overline{C_F(x)})$. 

These pairs are important for our model of the cube as a genotype space because they represent alternate strategies for viability. If the members of a pair are related in the partial order on strong components then only the strategy corresponding to the greater component is viable because the alleles in the lesser component require their complements for survival, which is not allowed. On the other hand, we will see that if there is no order relation between members of a nontrivial pair, then the pair splits the viable genotypes into disconnected clusters.

\begin{figure}
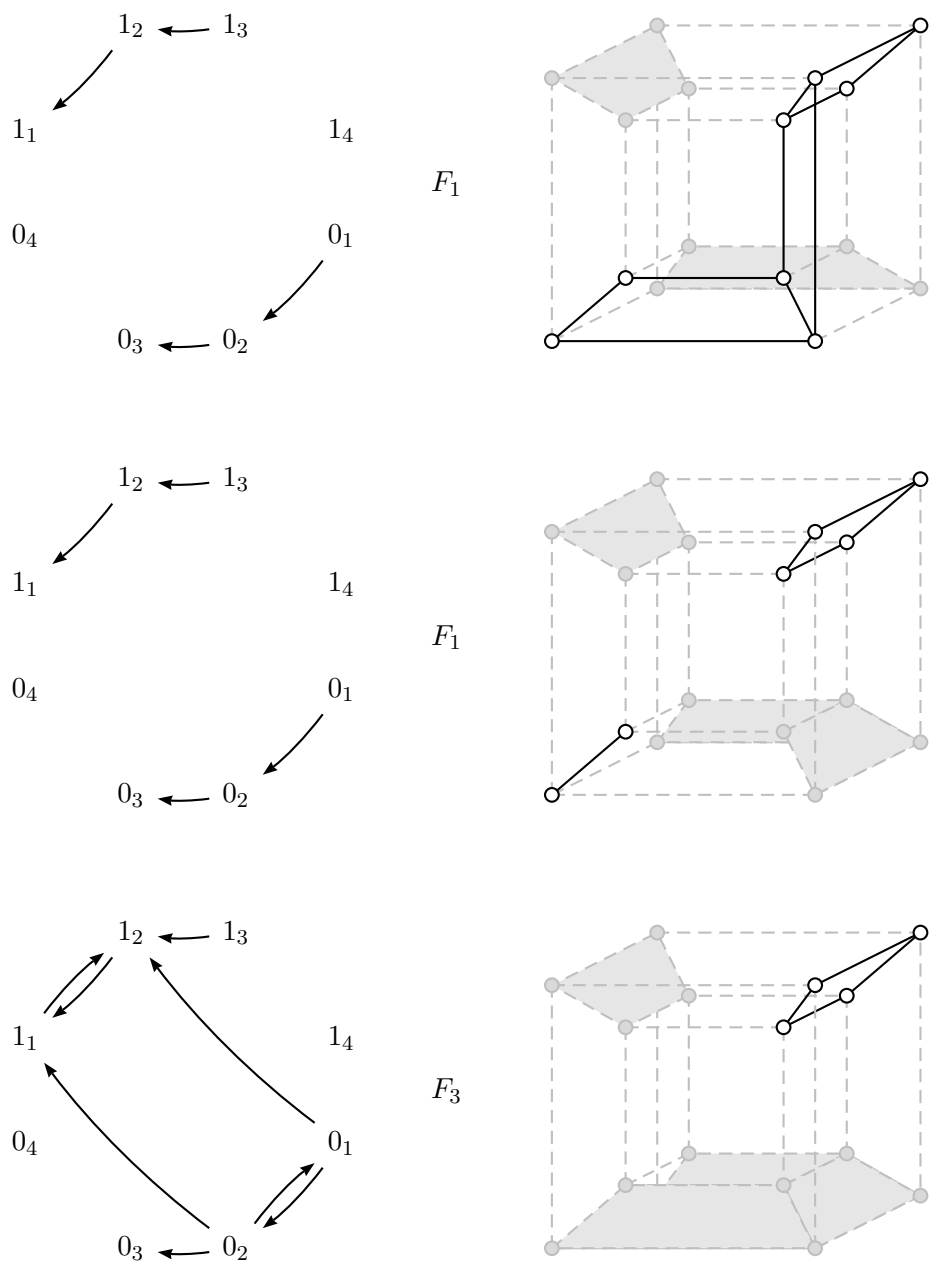
\label{fig}
\pspicture(-2,-3)(5,1.5)

\psset{unit=1.4,nodesep=5pt}
\rput(0,0){\rnode{1}{$1_1$}}
\rput(1,1){\rnode{2}{$1_2$}}
\rput(2,1){\rnode{3}{$1_3$}}
\rput(3,0){\rnode{4}{$1_4$}}
\rput(3,-1){\rnode{n1}{$0_1$}}
\rput(2,-2){\rnode{n2}{$0_2$}}
\rput(1,-2){\rnode{n3}{$0_3$}}
\rput(0,-1){\rnode{n4}{$0_4$}}
\ncarc[arrowsize=.1,arrowinset=.1]{->}{3}{2}
	\ncarc[arrowsize=.1,arrowinset=.1]{->}{n2}{n3}
\ncarc[arrowsize=.1,arrowinset=.1]{->}{n1}{n2}
	\ncarc[arrowsize=.1,arrowinset=.1]{->}{2}{1}
\rput(4,-.5){$F_1$}

\psset{origin={-5,2},unit=.5}
\psline(0,0)(5,0)
\psline[linestyle=dashed,linecolor=lightgray](0,0)(0,5)
\psline[linestyle=dashed,linecolor=lightgray](0,0)(2,1)(2,6)(7,6)(7,1)(2,1)
\pspolygon[fillstyle=solid,fillcolor=gray9,linestyle=dashed,linecolor=lightgray](2,1)(7,1)(5.6,1.8)(2.6,1.8)
\psline(0,0)(1.4,1.2)
\psline[linestyle=dashed,linecolor=lightgray](1.4,1.2)(1.4,4.2)
\psline[linestyle=dashed,linecolor=lightgray](1.4,4.2)(4.4,4.2)
\psline[linestyle=dashed,linecolor=lightgray](1.4,1.2)(2.6,1.8)
\psline(4.4,4.2)(4.4,1.2)(1.4,1.2)
\psline[linestyle=dashed,linecolor=lightgray](2.6,1.8)(2.6,4.8)(5.6,4.8)(5.6,1.8)
\psline(5,5)(7,6)
\psline(5,5)(4.4,4.2)
\psline[linestyle=dashed,linecolor=lightgray](5,0)(7,1)
\psline(5,0)(4.4,1.2)
\psline(4.4,4.2)(5.6,4.8)
\psline(5.6,4.8)(7,6)
\psline[linestyle=dashed,linecolor=lightgray](4.4,1.2)(5.6,1.8)
\pspolygon[fillstyle=solid,fillcolor=gray9,linestyle=dashed,linecolor=lightgray](0,5)(2,6)(2.6,4.8)(1.4,4.2)
\psline[linestyle=dashed,linecolor=lightgray](0,5)(5,5)
\psline(5,0)(5,5)
\pscircle[fillstyle=solid,fillcolor=white](0,0){.15}
\pscircle[fillstyle=solid,fillcolor=gray85,linestyle=solid,linecolor=lightgray](0,5){.15}
\pscircle[fillstyle=solid,fillcolor=white](5,5){.15}
\pscircle[fillstyle=solid,fillcolor=white](5,0){.15}
\pscircle[fillstyle=solid,fillcolor=gray85,linestyle=solid,linecolor=lightgray](2,1){.15}
\pscircle[fillstyle=solid,fillcolor=gray85,linestyle=solid,linecolor=lightgray](2,6){.15}
\pscircle[fillstyle=solid,fillcolor=white](7,6){.15}
\pscircle[fillstyle=solid,fillcolor=gray85,linestyle=solid,linecolor=lightgray](7,1){.15}
\pscircle[fillstyle=solid,fillcolor=white](1.4,1.2){.15}
\pscircle[fillstyle=solid,fillcolor=gray85,linestyle=solid,linecolor=lightgray](1.4,4.2){.15}
\pscircle[fillstyle=solid,fillcolor=white](4.4,4.2){.15}
\pscircle[fillstyle=solid,fillcolor=white](4.4,1.2){.15}
\pscircle[fillstyle=solid,fillcolor=gray85,linestyle=solid,linecolor=lightgray](2.6,1.8){.15}
\pscircle[fillstyle=solid,fillcolor=gray85,linestyle=solid,linecolor=lightgray](2.6,4.8){.15}
\pscircle[fillstyle=solid,fillcolor=white](5.6,4.8){.15}
\pscircle[fillstyle=solid,fillcolor=gray85,linestyle=solid,linecolor=lightgray](5.6,1.8){.15}

\endpspicture

\pspicture(-2,-3)(5,3)

\psset{unit=1.4,nodesep=5pt}
\psset{origin={0,-5}}
\rput(0,0){\rnode{1}{$1_1$}}
\rput(1,1){\rnode{2}{$1_2$}}
\rput(2,1){\rnode{3}{$1_3$}}
\rput(3,0){\rnode{4}{$1_4$}}
\rput(3,-1){\rnode{n1}{$0_1$}}
\rput(2,-2){\rnode{n2}{$0_2$}}
\rput(1,-2){\rnode{n3}{$0_3$}}
\rput(0,-1){\rnode{n4}{$0_4$}}
\ncarc[arrowsize=.1,arrowinset=.1]{->}{3}{2}
	\ncarc[arrowsize=.1,arrowinset=.1]{->}{n2}{n3}
\ncarc[arrowsize=.1,arrowinset=.1]{->}{n1}{n2}
	\ncarc[arrowsize=.1,arrowinset=.1]{->}{2}{1}
\rput(4,-.5){$F_1$}

\psset{origin={-5,2},unit=.5}
\psline[linestyle=dashed,linecolor=lightgray](0,0)(5,0)
\psline[linestyle=dashed,linecolor=lightgray](0,0)(0,5)
\psline[linestyle=dashed,linecolor=lightgray](0,0)(2,1)(2,6)(7,6)(7,1)(2,1)
\pspolygon[fillstyle=solid,fillcolor=gray9,linestyle=dashed,linecolor=lightgray](2,1)(7,1)(5.6,1.8)(2.6,1.8)
\psline(0,0)(1.4,1.2)
\psline[linestyle=dashed,linecolor=lightgray](1.4,1.2)(1.4,4.2)
\psline[linestyle=dashed,linecolor=lightgray](1.4,4.2)(4.4,4.2)
\psline[linestyle=dashed,linecolor=lightgray](1.4,1.2)(2.6,1.8)
\psline[linestyle=dashed,linecolor=lightgray](4.4,4.2)(4.4,1.2)(1.4,1.2)
\psline[linestyle=dashed,linecolor=lightgray](2.6,1.8)(2.6,4.8)(5.6,4.8)(5.6,1.8)
\psline(5,5)(7,6)
\psline(5,5)(4.4,4.2)
\psline[linestyle=dashed,linecolor=lightgray](5,0)(7,1)
\psline(4.4,4.2)(5.6,4.8)
\psline(5.6,4.8)(7,6)
\psline[linestyle=dashed,linecolor=lightgray](4.4,1.2)(5.6,1.8)
\pspolygon[fillstyle=solid,fillcolor=gray9,linestyle=dashed,linecolor=lightgray](5,0)(7,1)(5.6,1.8)(4.4,1.2)
\pspolygon[fillstyle=solid,fillcolor=gray9,linestyle=dashed,linecolor=lightgray](0,5)(2,6)(2.6,4.8)(1.4,4.2)
\psline[linestyle=dashed,linecolor=lightgray](0,5)(5,5)
\psline[linestyle=dashed,linecolor=lightgray](5,0)(5,5)
\pscircle[fillstyle=solid,fillcolor=white](0,0){.15}
\pscircle[fillstyle=solid,fillcolor=gray85,linestyle=solid,linecolor=lightgray](0,5){.15}
\pscircle[fillstyle=solid,fillcolor=white](5,5){.15}
\pscircle[fillstyle=solid,fillcolor=gray85,linestyle=solid,linecolor=lightgray](5,0){.15}
\pscircle[fillstyle=solid,fillcolor=gray85,linestyle=solid,linecolor=lightgray](2,1){.15}
\pscircle[fillstyle=solid,fillcolor=gray85,linestyle=solid,linecolor=lightgray](2,6){.15}
\pscircle[fillstyle=solid,fillcolor=white](7,6){.15}
\pscircle[fillstyle=solid,fillcolor=gray85,linestyle=solid,linecolor=lightgray](7,1){.15}
\pscircle[fillstyle=solid,fillcolor=white](1.4,1.2){.15}
\pscircle[fillstyle=solid,fillcolor=gray85,linestyle=solid,linecolor=lightgray](1.4,4.2){.15}
\pscircle[fillstyle=solid,fillcolor=white](4.4,4.2){.15}
\pscircle[fillstyle=solid,fillcolor=gray85,linestyle=solid,linecolor=lightgray](4.4,1.2){.15}
\pscircle[fillstyle=solid,fillcolor=gray85,linestyle=solid,linecolor=lightgray](2.6,1.8){.15}
\pscircle[fillstyle=solid,fillcolor=gray85,linestyle=solid,linecolor=lightgray](2.6,4.8){.15}
\pscircle[fillstyle=solid,fillcolor=white](5.6,4.8){.15}
\pscircle[fillstyle=solid,fillcolor=gray85,linestyle=solid,linecolor=lightgray](5.6,1.8){.15}

\endpspicture

\pspicture(-2,-4)(5,3)

\psset{unit=1.4,nodesep=5pt}
\rput(0,0){\rnode{1}{$1_1$}}
\rput(1,1){\rnode{2}{$1_2$}}
\rput(2,1){\rnode{3}{$1_3$}}
\rput(3,0){\rnode{4}{$1_4$}}
\rput(3,-1){\rnode{n1}{$0_1$}}
\rput(2,-2){\rnode{n2}{$0_2$}}
\rput(1,-2){\rnode{n3}{$0_3$}}
\rput(0,-1){\rnode{n4}{$0_4$}}
\ncarc[arrowsize=.1,arrowinset=.1]{->}{3}{2}
	\ncarc[arrowsize=.1,arrowinset=.1]{->}{n2}{n3}
\ncarc[arrowsize=.1,arrowinset=.1]{->}{1}{2}
	\ncarc[arrowsize=.1,arrowinset=.1]{->}{n2}{n1}
\ncarc[arrowsize=.1,arrowinset=.1]{->}{n1}{n2}
	\ncarc[arrowsize=.1,arrowinset=.1]{->}{2}{1}
\ncarc[arrowsize=.1,arrowinset=.1]{->}{n1}{2}
	\ncarc[arrowsize=.1,arrowinset=.1]{->}{n2}{1}
\rput(4,-.5){$F_3$}

\psset{origin={-5,2},unit=.5}
\psline[linestyle=dashed,linecolor=lightgray](0,0)(5,0)
\psline[linestyle=dashed,linecolor=lightgray](0,0)(0,5)
\psline[linestyle=dashed,linecolor=lightgray](0,0)(2,1)(2,6)(7,6)(7,1)(2,1)
\pspolygon[fillstyle=solid,fillcolor=gray9,linestyle=dashed,linecolor=lightgray](2,1)(7,1)(5.6,1.8)(2.6,1.8)
\psline[linestyle=dashed,linecolor=lightgray](1.4,1.2)(1.4,4.2)
\psline[linestyle=dashed,linecolor=lightgray](1.4,4.2)(4.4,4.2)
\psline[linestyle=dashed,linecolor=lightgray](1.4,1.2)(2.6,1.8)
\psline[linestyle=dashed,linecolor=lightgray](4.4,4.2)(4.4,1.2)(1.4,1.2)
\psline[linestyle=dashed,linecolor=lightgray](2.6,1.8)(2.6,4.8)(5.6,4.8)(5.6,1.8)
\psline(5,5)(7,6)
\psline(5,5)(4.4,4.2)
\psline[linestyle=dashed,linecolor=lightgray](5,0)(7,1)
\psline(4.4,4.2)(5.6,4.8)
\psline(5.6,4.8)(7,6)
\psline[linestyle=dashed,linecolor=lightgray](4.4,1.2)(5.6,1.8)
\pspolygon[fillstyle=solid,fillcolor=gray9,linestyle=dashed,linecolor=lightgray](5,0)(7,1)(5.6,1.8)(4.4,1.2)
\pspolygon[fillstyle=solid,fillcolor=gray9,linestyle=dashed,linecolor=lightgray](0,5)(2,6)(2.6,4.8)(1.4,4.2)
\pspolygon[fillstyle=solid,fillcolor=gray9,linestyle=dashed,linecolor=lightgray](0,0)(5,0)(4.4,1.2)(1.4,1.2)
\psline[linestyle=dashed,linecolor=lightgray](0,5)(5,5)
\psline[linestyle=dashed,linecolor=lightgray](5,0)(5,5)
\pscircle[fillstyle=solid,fillcolor=gray85,linestyle=solid,linecolor=lightgray](0,0){.15}
\pscircle[fillstyle=solid,fillcolor=gray85,linestyle=solid,linecolor=lightgray](0,5){.15}
\pscircle[fillstyle=solid,fillcolor=white](5,5){.15}
\pscircle[fillstyle=solid,fillcolor=gray85,linestyle=solid,linecolor=lightgray](5,0){.15}
\pscircle[fillstyle=solid,fillcolor=gray85,linestyle=solid,linecolor=lightgray](2,6){.15}
\pscircle[fillstyle=solid,fillcolor=white](7,6){.15}
\pscircle[fillstyle=solid,fillcolor=gray85,linestyle=solid,linecolor=lightgray](7,1){.15}
\pscircle[fillstyle=solid,fillcolor=gray85,linestyle=solid,linecolor=lightgray](1.4,1.2){.15}
\pscircle[fillstyle=solid,fillcolor=gray85,linestyle=solid,linecolor=lightgray](1.4,4.2){.15}
\pscircle[fillstyle=solid,fillcolor=white](4.4,4.2){.15}
\pscircle[fillstyle=solid,fillcolor=gray85,linestyle=solid,linecolor=lightgray](4.4,1.2){.15}
\pscircle[fillstyle=solid,fillcolor=gray85,linestyle=solid,linecolor=lightgray](2.6,1.8){.15}
\pscircle[fillstyle=solid,fillcolor=gray85,linestyle=solid,linecolor=lightgray](2.6,4.8){.15}
\pscircle[fillstyle=solid,fillcolor=white](5.6,4.8){.15}
\pscircle[fillstyle=solid,fillcolor=gray85,linestyle=solid,linecolor=lightgray](5.6,1.8){.15}

\endpspicture

\caption{An example for $n=4$ where the dimensions on the cube are ordered left/right, down/up, front/back, and in/out. The inviable (shaded) genotypes are left to give perspective and to demonstrate the way incompatibilities eliminate subcubes (shaded). 
Notice the progression from $F_1=\{ (0_2,1_3), (0_1, 1_2) \}$ where there are no cycles on strictly distinct alleles, to $F_2=F_1 \cup (1_1, 0_2) \}$ where there is a pair of complementary cycles, but no order between the cycles, to
$F_3=F_2 \cup (0_1, 0_2) \}$ where there is an order relation between the pair of cycles.}
\end{figure}

Consider the following example which is illustrated in figure \ref{fig}. We begin with a set of incompatibilities, $F_1=\{ (0_2,1_3), (0_1, 1_2) \}$. We see that all viable genotypes are connected in $Q_{F_1}$ and that $D_{F_1}$ has only trivial strong components. We then add an incompatibility to form $F_2=F_1 \cup (1_1, 0_2)$. In this case, there are two clusters of viable genotypes in $Q_{F_2}$ and a nontrivial strong component pair in $D_{F_2}$ that is unrelated in the partial order on strong components. By the addition of another incompatibility we get $F_3=F_2 \cup(0_1,0_2)$, where there is an order relation between the strong component pair; $C_{F_3}(0_1) \leq C_{F_3}(1_1)$ and every viable genotype in $Q_{F_3}$ has all the alleles $1_1$ and $1_2$ which comprise $C_{F_3}(1_1)$.

This motivates the following definition: a {\em splitting pair} is defined to be a pair of distinct, nontrivial strong components $(C,\overline{C})$ in $D_F$ that are unrelated in the order on strong components. Furthermore, we say that a set $S$ of loci contains the splitting pair $(C,\overline{C})$ if each allele $x\in C$ corresponds to a locus in $S$.

\begin{lem} \label{lem:conn}
Suppose $F$ is a satisfiable 2-formula and that $u$ and $v$ are any two viable genotypes in $\qf$, and let $S$ be the set of loci on which $u$ and $v$ differ. Then $u$ and $v$ are connected in $\qf$ if and only if $S$ contains no splitting pairs in $D_F$.
\end{lem} 

\begin{proof}
Suppose $S$ contains a splitting pair $(C,\overline{C})$ in $D_F$, on $k>1$ loci. Consider the subcubes of $Q^n$ which have alleles fixed to $C$ and $\overline{C}$ respectively. The Hamming distance between these subcubes is $k>1$ and there are no viable genotypes off these subcubes. Since $u$ and $v$ are in different subcubes, they are not connected by a viable mutational path and we see that splitting pairs are sensibly defined. 

On the other hand, suppose $S$ contains no splitting pairs in $D_F$. We must show that there is a path between $u$ and $v$ in $Q_F$. We think of each step toward $v$ as crossing out an element of $S$. Let $\overline{x}$ and $x$ be alleles of $u$ and $v$ respectively at an arbitrary locus in $S$. We construct a path starting from $u$ and ending in a genotype $w$ with allele $x$, such that each step decreases the distance to $v$, and such that each genotype on the path is viable. This is all that must be shown because each step crosses out an element of $S$ where $|S|<\infty$ and the locus of $x$ was arbitrary in $S$.

Suppose $x$ is an allele which has no outgoing edges in $D_F$. Then there are no incompatibilities involving $x$ and a mutation from a viable genotype with allele $\overline{x}$ to the neighboring genotype with allele $x$ will be a mutation between viable genotypes. More generally, the one step mutation from $\overline{x}$ to $x$ will not make a viable genotype inviable if and only if each of the $n-1$ shared alleles $y$ is such that 
$\overline{y} \notin L^+(x)$. 
Perhaps it is more intuitive to say that any allele in $L^+(x)$ other than $x$ is an allele that is shared by both genotypes. Since we are assuming that $v$ is viable with the allele $x$, we know that any allele $y \in L^+(x)$ is an allele of $v$. What must be shown is that we can follow a path of one step mutations between viable genotypes, which begins at $u$ and picks up all the alleles except $x$ in $L^+(x)$ which $u$ does not already contain. Then we will be able to make a one step mutation from $\overline{x}$ to $x$ that is between viable genotypes. We equate these steps in $Q_F$ with the removal of an allele from $L^+(x)$. Thus, we wish to remove all the alleles of $L^+(x)$.

Notice that the alleles of $L^+(x)$ are strictly distinct because $v$ is viable. Also, notice that if $u$ has any allele $y \in L^+(x)$, then $u$ has all the alleles of $L^+(y)$, else $u$ is inviable. Thus, all the alleles in the out-graphs of the alleles common to $u$ and $v$ are common to $u$ and $v$, so we need not make any of these mutations. Consequently, if we remove $D^+(y)$ from $D^+(x)$ for each allele $y$ common to $u$ and $v$, the remaining trimmed out subgraph of $D^+(x)$ will still be connected and any steps that we wish to make will concern only alleles in the trimmed out subgraph. Thus, we can assume that every allele in $L^+(x)$ that is common to $u$ and $v$ has already been removed. 

Suppose there is a nontrivial strong component $C \subset L^+(x)$. If there is an order relation between $C$ and $\overline{C}$, then $\overline{C} \leq C$ because $v$ is viable; but so too is $u$, so the alleles of $C$ are common to $u$ and $v$ and have already been removed from $L^+(x)$. 
Likewise, if there is no order relation between $C$ and $\overline{C}$, then $(C,\overline{C})$ is a splitting pair; so $u$ and $v$ share the alleles of $C$ by assumption. Consequently, we have shown that once the common alleles of $u$ and $v$ have been removed from $L^+(x)$, the trimmed out subgraph of $D^+(x)$ that remains is a directed tree. This means that only one step mutations are necessary to get from $u$ to a genotype with the alleles of $L^+(x)$.

To be precise, consider the leaves of the trimmed out subgraph. Any leaf $y$ is free of incompatibilities, so a mutation from a viable genotype with $\bar{y}$ to a genotype with $y$ will result in a viable genotype. Thus, by the viability of $u$, we see that the entire subcube of $Q^n$ formed by varying the loci of these leaves from $u$ is viable. Let $u^*$ be the genotype that results from making all these mutations. As we walk from $u$ to $u^*$ we trim the leaves and edges ending in those leaves from $D^+(x)$. The leaves of this new trimmed graph correspond to loci which can be varied from $u^*$ to form another viable subcube. Since $|L^+(x)|$ is finite, we can continue trimming the leaves until we are left with $x$ alone. At this point we will be free to follow the mutation from $\overline{x}$ to $x$ in $\qf$ and the path from $u$ to $w$ is complete. 

\end{proof}

Lemma~\ref{lem:conn} does not address the existence of viable genotypes
following the two possible strategies implied by each splitting pair, and
we turn to this issue next. As above, suppose $F$ is satisfiable and that $(C,\overline{C})$ is a splitting pair. Since $F$ is satisfiable, there is a viable genotype following at least one of the strategies $C$ or $\overline{C}$, say $u$ is following $\overline{C}$. We know that there is no $z \in C$ such that $z \lt \overline{z}$, else $(C,\overline{C})$ is not a splitting pair, so 
$ A = \bigcup_{z \in C} L^+(z) $
is a set of strictly distinct alleles and there are no incompatibilities among the pairs in $A$. Furthermore, none of the alleles of $u$ are incompatible. Thus, if there are not any incompatibilities $(x,y)$ such that $x$ is an allele of $A$ while $y$ is an allele of $u$ which is strictly distinct from the alleles of $A$ then the genotype with all the alleles of $A$ and alleles of $u$ assigned on the remaining loci will be viable. Suppose there were such an incompatibility. Then $x \lt \overline{y}$ so $\overline{y} \in A$, but then $y$ is not strictly distinct from the alleles of $A$ so we have a contradiction.

Thus, for every pair of potential strategies, there are viable genotypes with each stategy and we have the following 
\begin{cor}\label{cor:comps}
If $F$ is satisfiable and there are $k$ splitting pairs corresponding to $F$, then the number of clusters of viable genotypes in $Q_F$ is $2^k$. 
\end{cor}

\end{subsection}
\begin{subsection}{Proof of the theorem}

Let $\yn$ be the random variable that counts the number of splitting pairs in $D_F$. We will show that $\yn$ converges weakly to the Poisson random variable $\Psi$ with mean 
$$\la= - \frac{1}{2} \left( \ln{(1-c)}+c \right) .$$
We begin by discussing $i$-cycles. Let $\Gin$ to be the set of complementary $i$-cycle pairs in $D_F$ and define
$$\Gn=\bigcup_{i=2}^n \Gin.$$ 
For each $\al \in \Gn$ let 
$$I_\al=\left\{\begin{array}{ll}
1 & \text{if } \al \in D_F \\
0 & \text{else. }
\end{array}\right.$$
Then 
$$\xin=\sum_{\al \in \Gin} I_\al$$
is the random variable that counts the number of $i$-cycle pairs in $D_F$ and  
$$\xn=\sum_{i=2}^n \xin$$
counts the total number of simple cycle pairs in $D_F$. 

Notice that the number of directed $i$-cycles on strictly distinct alleles is $$\binom{n}{i} (i-1)! 2^i = \frac{(n)_i}{i} 2^i$$
where $(n)_i=n(n-1)\ldots (n-(i-1))$. For each such cycle 
$$\alpha= y_1\to y_2 \to \ldots \to y_i \to y_1$$
that is present in $D_F$, the complementary cycle 
$$\bar{\alpha}= \bar{y}_1 \leftarrow \bar{y}_2 \leftarrow \ldots \leftarrow \bar{y}_i \leftarrow \bar{y}_1$$ 
is also in $D_F$. Thus, the expected number of $i$-cycle pairs is
$$\muin = E \left( \xin \right) = \sum_{\alpha\in\Gin} p^i 
= \frac{(n)_i}{i} 2^{i-1} p^i \lton \frac{c^i}{2i}.$$
If $c<1$, then we also have that the expected total number of simple cycles is 
$$\la_n = E \left( \xn \right) =  \sum_{i=2}^n \muin
\lton - \frac{1}{2} \left( \ln{(1-c)} + c \right).$$

This immediately implies that if $c<1$ and $\omega=\omega(n) \to \infty$, then
\begin{align}
\label{mtail}
P \left( \sum_{i=\omega}^n \xin > 0 \right) 
\leq E\left( \sum_{i=\omega}^n \xin \right) \lton 0,
\end{align}
which suggests that truncated random variables should be enough for weak convergence. Thus, we define
$$\tmn = \sum_{i=2}^m X_i .$$

One way to prove the weak convergence of $\yn$ to $\Psi$ is to show that if $A\subset Z^+$ is any subset of the positive integers and $\ep>0$ is any constant, then there is an $n$ such that
\begin{align*}
\left| P(\yn \in A) - P(\Psi \in A) \right| < \ep .
\end{align*}

Let $\la_{m,n}$ be the expectation of $\tmn$ and $\Psi_{m,n}$ be the Poisson random variable with mean $\la_{m,n}$. Then the difference above can be bounded above by
\begin{align}
& \left| P(\yn \in A) - P(\tmn \in A) \right| \\
&+ \left| P(\tmn \in A) - P(\Psi_{m,n} \in A) \right| \\
&+ \left| P(\Psi_{m,n} \in A) - P(\Psi \in A) \right| .
\end{align}
Thus, to prove the weak convergence of $\yn$ to $\Psi$ it is enough to show that for every $\ep>0$ there is an $m$ and an $n$ which may depend on $m$, such that (2), (3) and (4) are each less than $\ep/3$.

We can certainly find $m$ and $n$ for (4) so we focus on (2) and (3).
Notice that (2) is bounded above by $P(\yn \neq \tmn)$, which is in turn bounded above by the expected number of events which distinguish $\yn$ from $\tmn$. To understand these distinguishing events, notice that if there are no $i$-cycles for $i>m$, no compound cycles of any length, and no paths between any complementary sets of alleles of size at most $m$, then  $\yn = \tmn$. Thus, the lemma below provides the desired bound for (2). 

\begin{lem}\label{lem:exps}
Let $p=c/(2n)$ and $c<1$. Also, let $m<\infty$ be fixed and 
$\ep_m = \sum_{i=m+1}^\infty \frac{c^i}{2i}.$ 
Then 
\begin{itemize}
\item[(i)] The expected number of complementary simple cycle pairs in $D_F$ of length greater than $m$ is less than $\ep_m$.
\item[(ii)] The expected number of compound cycles of $D_F$ on strictly distinct alleles is $O(n^{-1})$. 
\item[(iii)] The expected number of complementary simple cycle pairs in $D_F$ that are comparable in the order relation on strong components is $O(n^{-1})$.
\end{itemize}
\end{lem}

\begin{proof}
Item $(i)$ follows from (\ref{mtail}) above.

Item $(ii)$ follows from the observation that any compound cycle contains at least a simple cycle $\alpha$ and a path $\beta$ with first and last alleles in $\alpha$. 
To see this, suppose $\alpha \in \Gin$ and $\beta$ is a path with $j>0$ edges and both endpoints in $\al$. Then the expected number of such configurations in $D_F$ is bounded above by
$$\binom{n}{i}2^{i- 1}(i-1)!\binom{n-i}{j-1}2^{j-1}(j-1)! i^2p^{i+j} 
< i(2n)^{i+j-1}p^{i+j} = c^{i+j}i/(2n).$$
Summing over $1\leq j \leq n-i$ and then over $2\leq i\leq n$ gives the desired result. 

For $(iii)$, notice that the expected number of complementary simple cycle pairs that are comparable in the order relation on strong components is exactly the same as the expected number of cycle and path pairs described above. Indeed, in either case, there are $i(i-1)$ paths of length one and $i^2$ choices for the endpoints of longer paths and exactly the same alleles to choose from for intermediary alleles on the longer paths. 
\end{proof}

The supremum of (3)
is half the total variation between the law of $\tmn$ and $\Psi_{m,n}$. Thus, the bound for (3) follows from
\begin{lem}\label{lem:trunc}
Let $m<\infty$ be fixed. Then the total variation between the law of $\tmn$ and Poisson random variable with mean $\la_{m,n}$ is $O(n^{-1})$. 
\end{lem}
\begin{proof}
Let $\Gamma_\al$ be the set of all complementary cycle pairs that are distinct from $\al$ and yet not independent of $\al$. Also, let $E(I_\al)=p_\al$. The local Chen-Stein method~\cite{BHJ} gives the following upper bound on the total variation between the law of $\tmn$, and $\Psi_{m,n}$: 
\begin{align*}
\min \{1, \lanm^{-1} \} 
\left(
\sum_{ \al \in \Gn } p_\al^2 + 
\sum_{\al \in \Gn,\be \in \Gamma_\al} 
(p_\al p_\be + E(I_\al I_\be)) 
\right).
\end{align*}

Since $\lanm$ is bounded below we need only show that each of the sums above is $O(n^{-1})$. Notice also that, for each fixed $2\leq i\leq n$, $\muin$ is bounded below in $n$ and that the sums above strictly contain the related sum for any $\xin$. Thus, if the sum above is $O(n^{-1})$, then so too is the corresponding sum for each $\xin$. Consequently, the fact that $\xin$ converges weakly to the 
Poisson random variable with mean $c^i/(2i)$ follows from the proof given here. 

Consider the first sum
\begin{align*}
\sum_{ \al \in \Gn } p_\al^2 =
\sum_{i=2}^n \binom{n}{i}2^{i-1}(i-1)! p^{2i} 
< \sum_{i=2}^\infty \frac{1}{n^i} .
\end{align*}
Since this sum is $O(n^{-2})$, we can address the second sum. If $I_\al$ and $I_\be$ are not independent, then they share some edge and the event that $I_\al=1$ increases the probability that $I_\be=1$, so 
$$E(I_\al I_\be) > p_\al p_\be$$ 
and we need only show convergence of the sum 
\begin{align}\label{depSum}
\sum_{\al \in \Gn,\be \in \Gamma_\al} E(I_\al I_\be). 
\end{align}

Let $\al \in \Gamma_i$ and let $\be$ be such that the edges of $\be$ which are shared with those of $\al$ are found on $j$ disjoint paths in each of its members and such that there are $k$ alleles on each member of $\be$ which are not on the paths common to both $\al$ and $\be$. Then, the probability that any such fixed pair $(\al,\be)$ is present in $D_F$ is
$p^{i+j+k}.$

We will count the number of such pairs by first examining the number of possible shared paths. Notice first that $1 \leq j \leq \lfloor i/2 \rfloor$ and by choosing $2j$ alleles in $\al$ we count the number of possible endpoints of the shared paths in one of the cycles of $\al$. Between any two endpoints that have no other endpoints between them, the edges connecting them are either a shared path or not, but once this choice is made for one such pair, it is made for all the endpoints. Thus, the number of possible shared paths is
$$2\binom{i}{2j} 
< 2 \frac{ (i)_{2j} }{ 2^j j! 2^{j-1} (j-1)! } 
< 4 \frac{(i)_j}{j!} j \frac{i^j}{4^j j!}.$$
The $k$ alleles on each member of $\be$ could be on $\al$ or not and the ordering of these alleles is interspersed with the the ordering of the shared paths. There are less than $(2n)^k$ orderings for the $k$ alleles and at most $j!$ orderings allowed for the shared paths. Once these have been determined, we must choose where the $j$ ordered paths are placed among the $k$ ordered alleles. If there were no restrictions on interspersal, then this would be equivalent to forming weak compositions of $k$ into $j$ parts. Thus, we get the following upper bound on the number of ways to finish $\be$ once $\al$ and the shared paths have been determined:
$$(2n)^k j! \binom{ k+j-1 }{ j-1 }.$$
Hence, for fixed $i$, $j$ and $k$, an upper bound on the expected number of such ordered pairs $(\al,\be)$ in $D_F$ is
\begin{align*}
&
4 \frac{ (n)_i }{ i } 2^i
j \binom{i}{j} \left( \frac{i}{4} \right)^j 
(2n)^k \binom{ k+j-1 }{ j-1 } p^{i+j+k}
\\ & =
4 \frac{ (n)_i }{ i } \left( \frac{c}{n} \right)^i  
j \binom{i}{j} \left( \frac{ci}{8n} \right)^j 
\binom{ k+j-1 }{ k } c^k.
\end{align*}
The sum over the terms involving $k$ can be bounded above by the sum over all $k \geq 0$:
\begin{align*}
\sum_{k=0}^\infty \binom{ k+j-1 }{ k } c^k = \left( \frac{1}{1-c} \right)^j .
\end{align*}
The sum over $j$ is bounded as follows
\begin{align*}
\sum_{j=1}^{\lfloor i/2 \rfloor}
j \binom{i}{j} \left( \frac{ci}{8(1-c)n} \right)^j 
&=\sum_{j=1}^{\lfloor i/2 \rfloor}
i \binom{i-1}{j-1} 
\left( \frac{ci}{8(1-c)n} \right)^j 
\\ & \leq
\left( \frac{ci^2}{8(1-c)n} \right) 
\sum_{j=0}^i
\binom{i}{j} 
\left( \frac{ci}{8(1-c)n} \right)^j 
\\ & =
p\frac{i^2}{4(1-c)} \left( 1+\frac{ci}{8(1-c)n} \right)^i.
\end{align*}
Using the bounds 
\begin{align*}
1+x \leq e^x
\quad\text{ and }\quad
\frac{ (n)_i }{ n^i} 
=\prod_{l=1}^{i-1} \left( 1- \frac{l}{n} \right)
\leq e^{ -\frac{i^2}{2n} }
\end{align*}
we can get an upper bound for the sum over $2\leq i\leq m$
\begin{align}\label{sumi}
\frac{p}{1-c} \sum_{i=2}^m  
i c^i \exp{\left[ \left( -1  + \frac{c}{8(1-c)} \right) \frac{i^2}{2n} \right] }.
\end{align}
Since $m$ is fixed, this sum is $O(n^{-1})$, which completes the proof of the lemma.
\end{proof}

This also proves that $Y_n$ converges weakly to $\Psi$ and consequently that the probability that there is no splitting pair converges to $e^{-\la}$. 
Moreover, since the formula $F_p^n$ is satisfiable a.a.s., $N_n$ and $2^{Y_n}$ both converge weakly to $N=2^\Psi$. Thus, we have completed the proof of the theorem.

\end{subsection}
\begin{subsection}{Discussion}

We find the structure of this fitness landscape to be interesting because, while the number of viable genotypes grows exponentially with the number of loci, the number of clusters of these genotypes is relatively small and may be finite in expectation. Computer simulations have been run~\cite{GPG} that suggest this is the case, but we have not been able to determine whether or not the expectation of $N_n$ converges. This appears to present a challenge because one must consider correlation between cycles in $D_F$ of unbounded length. For a discussion of correlation between undirected cycles of bounded length, see~\cite{Bol}; or for a result concerning directed cycles, see~\cite{Pal}. 

\newpage

\end{subsection}
\end{section}

\end{document}